\documentclass[a4paper]{llncs}

\usepackage[utf8]{inputenc} 
\usepackage[T1]{fontenc}    
\usepackage{hyperref}       
\usepackage{url}            
\usepackage{booktabs}       
\usepackage{amsfonts}       
\usepackage{nicefrac}       
\usepackage{microtype}      

\usepackage{amssymb}
\usepackage{amsmath}
\usepackage{graphicx}

\usepackage{graphicx}
\usepackage{makecell}
\usepackage{float}
\newcommand{\tensor}[1]{\boldsymbol{\mathcal{#1}}}

\newcommand{\mat}[1]{\mathbf{#1}}
\newcommand{\vect}[1]{\mathbf{#1}}

\renewcommand{\figurename}{Fig.} 
\renewcommand{\tablename}{Tab.}

\begin{document}
\mainmatter 
\title{Rank Minimization on Tensor Ring: A New Paradigm in Scalable Tensor Decomposition and Completion}
\author{\small Longhao Yuan$^{1,2}$ \and Chao Li$^2$ \and Danilo Mandic $^{5}$ \and Jianting Cao$^{4,1}$\and Qibin Zhao$^{2,3}$}
\institute{ Graduate School of Engineering, Saitama Institute of Technology, Japan
\and Tensor Learning Unit, RIKEN Center for Advanced Intelligence Project (AIP), Japan
\and School of Automation, Guangdong University of Technology, China
\and School of Computer Science and Technology, Hangzhou Dianzi University, China
\and Department of Electrical and Electronic Engineering, Imperial College London, United Kingdom
\\ \url{ {longhao.yuan, chao.li, qibin.zhao}@riken.jp,
 cao@sit.ac.jp, d.mandic@imperial.ac.uk}}
\maketitle

\begin{abstract}
In low-rank tensor completion tasks, due to the underlying multiple large-scale singular value decomposition (SVD) operations and rank selection problem of the traditional methods, they suffer from high computational cost and high sensitivity of model complexity. In this paper, taking advantages of high compressibility of the recently proposed tensor ring (TR) decomposition, we propose a new model for tensor completion problem. This is achieved through introducing convex surrogates of tensor low-rank assumption on latent tensor ring factors, which makes it possible for the Schatten norm regularization based models to be solved at much smaller scale. We propose two algorithms which apply different structured Schatten norms on tensor ring factors respectively. By the alternating direction method of multipliers (ADMM) scheme, the tensor ring factors and the predicted tensor can be optimized simultaneously. The experiments on synthetic data and real-world data show the high performance and efficiency of the proposed approach.
\end{abstract}

\section{Introduction}
Tensor decomposition aims to find the latent factors of tensor valued data (i.e. the generalization of multi-dimensional arrays), thereby casting large-scale tensors into a multilinear tensor space of low-dimensionality (very few degree of freedom designated by the rank). Tensor factors can then be considered as latent features of data, and in this way can represent the data economically and predict missing entries when the data is incomplete. The specific form and operations among latent factors defines the type of tensor decomposition. A variety of tensor decomposition models have been applied in diverse fields such as machine learning \cite{novikov2015tensorizing,anandkumar2014tensor,kanagawa2016gaussian} and signal processing \cite{zhou2016linked,cong2015tensor}. Tucker decomposition and CANDECOMP/PARAFAC (CP) decomposition are classical tensor decomposition models, which have been studied for nearly half a century \cite{kolda2009tensor,tucker1966some,harshman1970foundations}. 

In recent years, the concept of tensor networks has been proposed and has become a powerful and promising aspect of tensor methodology \cite{cichocki2016tensor,cichocki2017tensor}. One of the most recent and popular tensor networks, named the matrix product state/tensor-train (MPS/TT), is studied across disciplines owing to its super compression and computational efficiency properties \cite{oseledets2011tensor,novikov2015tensorizing}. For a tensor of $N$-dimensions, the most significant property of TT decomposition is that the space complexity will not grow exponentially with $N$, thus providing a natural remedy for the `curse of dimensionality', while the number of parameters of Tucker decomposition is exponential in $N$. Although the CP decomposition is a highly compact representation which has desirable property of being linear in $N$, it has difficulties in finding the optimal latent tensor factors. To address these issues, recent studies propose a generalization of TT decomposition, termed the tensor ring (TR) decomposition, in order to relax the rank constraint of TT, thus offering an enhanced representation ability, latent factors permutation flexibility (i.e. tensor permutation is directly related to the permutation of tensor factors) and structure information interpretability (i.e. each tensor factor can represent a specific feature of original tensor) \cite{zhao2016tensor,zhao2018learning}.

Tensor completion aims to recover an incomplete tensor from partially observed entries. The theoretical lynchpin in matrix or tensor completion problems is the low-rank assumption, and tensor completion has been applied in various applications such as image/video completion \cite{liu2013tensor,zhao2015bayesian}, recommendation systems \cite{karatzoglou2010multiverse}, link prediction \cite{ermics2015link}, compressed sensing \cite{gandy2011tensor}, to name but a few. Since the determination of tensor rank is an NP-hard problem\cite{hillar2013most,kolda2009tensor}, many tensor low-rank surrogates were proposed for tensor completion. One such surrogate is the Schatten norm (a.k.a. nuclear norm, or trace norm),  which is defined as the sum of singular values of a matrix, and is the most popular convex surrogate for rank regularization. Unlike matrix completion problems, the Schatten norm model of a tensor is hard to formulate. Recent studies mainly focus on two convex relaxation models of tensor Schatten norm, the `overlapped' model \cite{liu2013tensor,tomioka2013convex,cheng2016scalable,signoretto2014learning,imaizumi2017tensor} and the `latent' \cite{tomioka2013convex,guo2017efficient} model. 

The work in \cite{tomioka2013convex} first proposes the `latent' norm model and shows that the mean square error of a `latent' norm method scales no greater than the `overlapped' norm method. Under the low-rank regularization of the latent model, the tensor does not need to be low-rank at every mode, which is considered as a more flexible constraint. Both models do not need to specify the rank of decompositions, and the rank of tensor is optimized to be minimal subject to the equivalence of observed elements. However, the two methods need to perform multiple SVD operations on the matricization of tensors, and the computational complexity grows exponentially with tensor dimension. Other tensor completion algorithms, like alternating least squares (ALS) \cite{grasedyck2015variants,wang2017efficient} and gradient-based algorithms \cite{yuan2017completion,acar2011scalable}, need to specify the rank of the decompositions beforehand, which leads to annoyed parameter tuning problems. In addition, the completion performance of tensor completion algorithms is mainly affected by rank selection, the number of observed entries and tensor dimensions.

In this paper, in order to tackle the high computational cost and the sensitivity to rank selection problems that most proposed algorithms experience, we propose a new tensor completion model based on the tensor ring decomposition. Our main contributions are listed below:
\begin{itemize}
\item The relation between low-rank assumption on tensor and latent factors is theoretically explained , and the low-rank surrogate on latent factors of tensor ring decomposition is introduced.

\item We formulate the TR overlapped low-rank factor (TR-OLRF) model and the TR latent low-rank factor (TR-LLRF) model, then the two models are solved efficiently by the ADMM algorithm.

\item We conduct several experiments and obtain the high performance and high efficiency by using our algorithms. In addition, the experiments results also show that our algorithms are robust to rank selection and data dimensionality.
\end{itemize}

\section{Preliminaries}

\subsection{tensor ring decomposition}
\begin{figure}[htb]
\begin{center}
\includegraphics[width=0.5\linewidth]{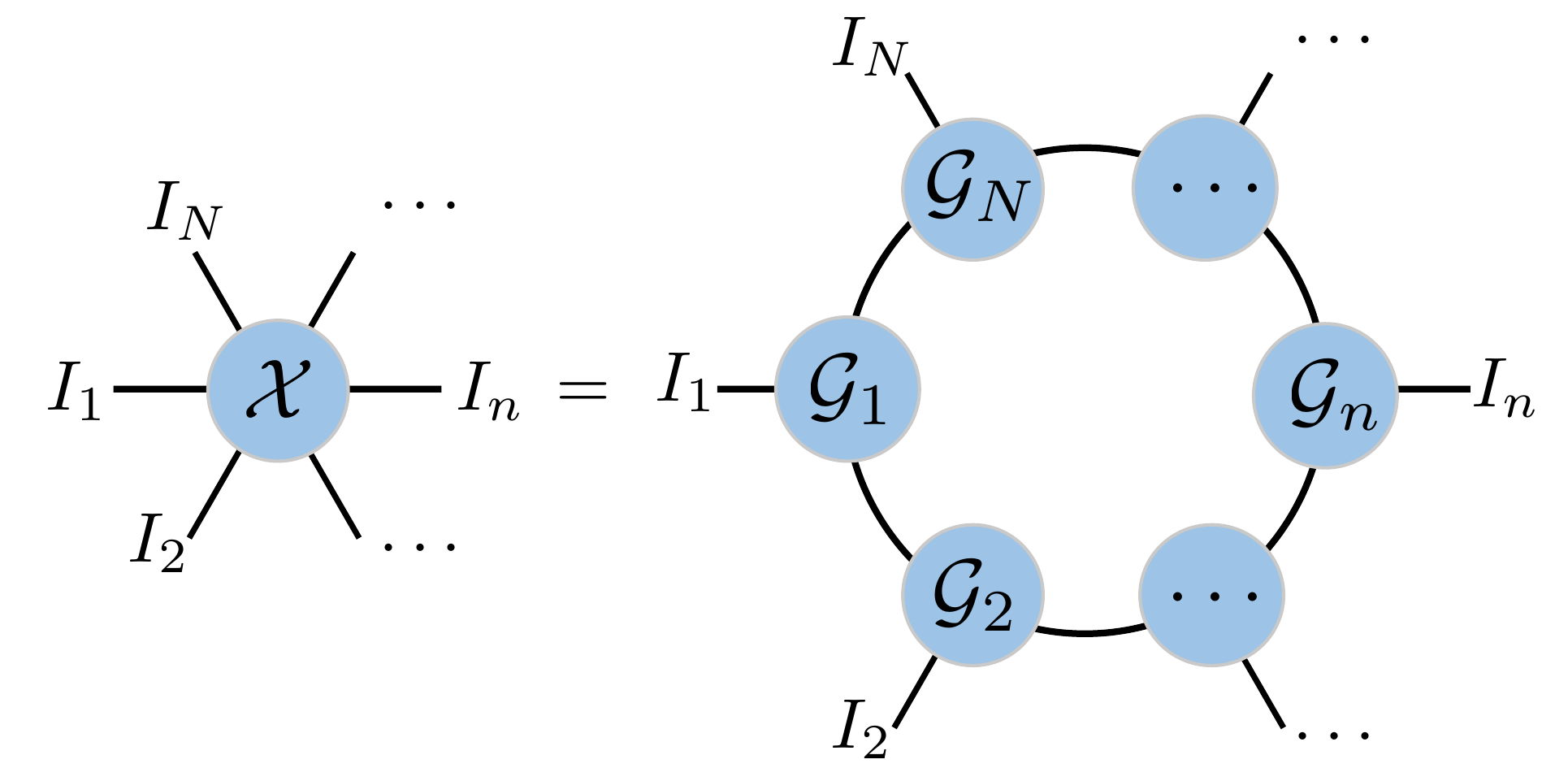}
\end{center}
\caption{TR decomposition}\label{tr_d}
\end{figure}
tensor ring (TR) decomposition is a more general decomposition than tensor-train (TT) decomposition, and it represents a tensor with large dimension by circular multilinear products over a sequence of low dimension cores. All of the cores corresponding to TR decomposition are order-three tensors, and are denoted by $\tensor{G}_{n} \in\mathbb{R}^{R_{n} \times I_{n} \times R_{n+1}}$, $n=1,\ldots,N$. The decomposition diagram is shown in Fig. \ref{tr_d}. In the same way as TT, the TR decomposition linearly scales to the dimension of the tensor, thus it can overcome the `curse of dimensionality'. For simplicity, we define $\{ \tensor{G}_n\}_{n=1}^N:=\{\tensor{G}_{1},\tensor{G}_{2},\ldots,\tensor{G}_{N}\}$ to represent a set of tensor cores. The syntax $\{R_n\}_{n=1}^N$ denotes TR-rank which controls the model complexity of TR decomposition. The TR decomposition relaxes the rank constraint on the first and last core of TT to $R_1=R_{N+1}$, while the original constraint on TT is rather stringent, i.e., $R_1=R_{N+1}=1$. TR applies trace operation and all the core tensors are constrained to be third-order equivalently. In this case, TR can be considered as a linear combination of TT and thus offers a more powerful and generalized representation ability than TT. The element-wise relation and global relation of TR decomposition and the original tensor is given by equations (\ref{tr_relation1}) and (\ref{tr_relation2}):
{\setlength\abovedisplayskip{5pt}
\setlength\belowdisplayskip{-10pt}
\begin{equation}
\label{tr_relation1}
\tensor{X}(i_1,i_2,\ldots,i_N)=\text{Trace}\left \{  \prod_{n=1}^N \mat{G}_n(i_n) \right \},
\end{equation}}

\begin{equation}
\label{tr_relation2}
\Delta_n(\mat{X})=\Gamma_2(\mat{G}_n)\Delta_2(\mat{G}_{\neq n})^T,
\end{equation}
where $\text{Trace}\{ \cdot \}$ is the matrix trace operator, $ \mat{G}_n(i_n)\in\mathbb{R}^{R_n\times R_{n+1}}$ is the $i_n$th mode-$2$ slice matrix of $\tensor{G}_n$, which also can be denoted by $\tensor{G}_n(:,i_n,:)$. $\tensor{G}_{\neq n}\in\mathbb{R}^{R_{n+1}\times \prod_{i=1, i\neq n}^N I_i \times R_n}$ is a subchain tensor by merging all cores except the $n$th core tensor, i.e., $\tensor{G}_{n+1},\ldots,\tensor{G}_N, \tensor{G}_1, \ldots ,\tensor{G}_{n-1}$. $\Gamma_n(\cdot)$ is the mode-$n$ matricization operator of a tensor, i.e., if $\tensor{X}\in\mathbb{R}^{I_1\times I_2\times\cdots \times I_N}$, then $\Gamma_n(\mat{X})\in\mathbb{R}^{I_n \times  {I_1 \cdots I_{n-1} I_{n+1} \cdots I_N}}$. $\Delta_n(\cdot)$ is another type of mode-$n$ matricization operator of a tensor, e.g., if $\tensor{X}\in\mathbb{R}^{I_1 \times  \cdots \times I_N}$, then $\Delta_n(\mat{X})\in\mathbb{R}^{I_n \times  {I_{n+1} \cdots I_{N} I_{1} \cdots I_{n-1}}}$.

\subsection{Tensor completion by Schatten norm regularization}

The low-rank tensor completion problem can be formulated as:
\begin{equation}
\label{bm}
\min \limits_{\tensor{X}} \ \  \text{Rank}(\tensor{X}),\ s.t. \ P_{\Omega}(\tensor{X})=P_{\Omega}(\tensor{T}),
\end{equation}
and the model can be written in a unconstrained form by:
\begin{equation}
\label{bmc}
\min \limits_{\tensor{X}} \ \  \text{Rank}(\tensor{X})+\frac{1}{\lambda}\Vert P_{\Omega}(\tensor{X})-P_{\Omega}(\tensor{T})\Vert_F^2,
\end{equation}
where $\tensor{X}$ is the low-rank approximation tensor, $\text{Rank}(\cdot)$ is a rank regularizer, $P_{\Omega}(\tensor{T})$ denotes all the observed entries w.r.t. the set of indices of observed entries represented by $\Omega$, and $\Vert \cdot \Vert_F$ is the Frobenius norm.
For the low-rank tensor completion problem, determining the rank of a tensor is an NP-hard problem. Work in \cite{liu2013tensor} and \cite{signoretto2014learning} extends the concept of low-rank matrix completion and defines the tensor rank as the sum of the rank of mode-$n$ matricization of the tensor. This surrogate is named `overlapped' model, and it simultaneously regularizes all the mode-$n$ matricizations of a tensor into low-rankness by Schatten norm. In this way, we can define the rank of a tensor as:
\begin{equation}
\label{overlapped}
\text{Rank}(\tensor{X}):=\sum\nolimits_{n=1}^N\Vert \Gamma_n(\mat{X}) \Vert_*, 
\end{equation}
where $\Vert \cdot \Vert_*$ denotes the Schatten norm. 

Another surrogate of tensor rank, named `latent' low-rank, has been proposed and studied recently. In \cite{tomioka2013convex}, the `latent' model considers the original tensor as a summation of several latent tensors and assumes that each latent tensor is low-rank in a specific mode:
\begin{equation}
\label{latent}
\text{Rank}(\tensor{X}):=\sum\nolimits_{n=1}^N \Vert \Gamma_n(\mat{W}_n) \Vert_*, \ s.t.\ \tensor{X}=\sum\nolimits_{n=1}^N \tensor{W}_n.
\end{equation}

This convex surrogate is more flexible as it can fit the tensor well if the tensor does not have low-rankness in all modes. The completion algorithms based on these two models are shown to have fast convergence and good performance when data size is small. However, when we need to deal with large-scale data, the multiple SVD operations will be intractable due to high computational cost.

\subsection{Tensor completion by tensor decomposition}
Some other existing tensor completion algorithms do not employ a low-rank constraint to the tensor, and thus they do not find the low-rank tensor directly, instead, they try to find the low-rank representation (i.e. tensor factors) of the incomplete data by observed entries, then the obtained latent factors are used to predict the missing entries. The completion problem is set as a weighted least squares model, e.g., the tensor completion model based on TR decomposition is formulated below:
\begin{equation}
\label{tdc}
\min\limits_{\{\tensor{G}_n\}} \Vert \tensor{W}*(\tensor{T}-\tensor{X}(\tensor{G}_1,\ldots,\tensor{G}_N) \Vert_F^2,
\end{equation}
where $*$ is the Hadamard product of two tensors of same size, $\tensor{X}(\tensor{G}_1,\ldots,\tensor{G}_N)$ is the tensor generated by the tensor factors. $\tensor{W}$ is a weight tensor which is the same size as $\tensor{X}$, it records the indices of the observed entries of $\tensor{T}$, and every entry of $\tensor{W}$ satisfies $P_{\Omega}(\tensor{W})=1$ and $P_{\bar{\Omega}}(\tensor{W})=0$.

Based on solving tensor factors of different tensor decompositions, many tensor completion algorithms have been proposed, e.g., weighted CP \cite{acar2011scalable}, weighted Tucker \cite{filipovic2015tucker}, and weighted TT \cite{yuan2017completion}, TR-ALS \cite{wang2017efficient}. However, usually these algorithms are solved by gradient-based method or alternating least squares method, they are shown to suffer from low convergence speed and high computational cost. In addition, the performance of these methods is sensitive to rank selection.

In this paper, we make virtue of applying both `overlapped' approach and `latent' approach of structured Schatten norms, and aim to formulate a new tensor completion model. The main idea is to give a low-rank constraint on latent factors of a tensor. In this way, we only need to calculate SVD on the tensor factors instead of the whole scale of data. At the same time, low-rankness constraint on tensor factors will regularize the tensor factors to low-rank, and in doing so it will solve the problem of rank selection. The next section we presents our proposed method based on both `overlapped' and `latent' tensor low-rank models.

\section{Low-rankness on tensor factors}

We propose a new definition on low-rank tensor, which gives the low-rankness on the decomposition factors of a tensor, for TR decomposition, the low-rank model is formulated as:
\begin{equation}
\label{mod}
\begin{aligned}
\min \limits_{\{\tensor{G}_n\}_{n=1}^N,\tensor{X}}  \ \sum_{n=1}^N \text{Rank}(\tensor{G}_n) + \frac{\lambda}{2}\Vert \tensor{X}- \tensor{Z}(\{\tensor{G}_n \}_{n=1}^N)\Vert_F^2 , \ s.t. \ P_\Omega(\tensor{X})=P_\Omega(\tensor{T}),
\end{aligned}
\end{equation}
where $\tensor{Z}(\{\tensor{G}_n\}_{n=1}^N)$ denotes the tensor approximated by core tensors $ \{\tensor{G}_n\}_{n=1}^N$. We formulate the low-rank assumption of the core tensors by equation (\ref{overlapped}) and (\ref{latent}).

We need firstly to deduce the relation of tensor rank and tensor factor rank, which can be explained by the below theorem:

$Theorem\ 1$: For $n=1,\ldots,N$, $\sum_{i=1}^3 \text{Rank}(\Gamma_i(\mat{G}_n))) \geq \text{Rank}(\Delta_n(\mat{X}))$ .

$Proof$: For $n=1,\ldots,N$, from equation (\ref{tr_relation2}), we can infer $\text{Rank}(\Delta_n(\mat{X}))) \leq \text{Rank}(\Gamma_2(\mat{G}_n))) \leq \sum_{i=1}^3 \text{Rank}(\Gamma_i(\mat{G}_n))$.

The above theorem proves the relation between the ranks of tensor $\tensor{X}$ and core tensors $\{ \tensor{G}_n\}_{n=1}^N$. Since $\sum_{i=1}^3 \text{Rank}(\Gamma_i(\mat{G}_n)))$ is an upper bound of the mode-$n$ matricization of tensor $\tensor{X}$, we can take assumption that $\tensor{G}_n$ has a low-rank structure. This can largely decrease the computational complexity compared to other algorithms which give low-rank assumption on overlapped tensors or latent tensors. In a similar way, we can deduce that the sum of latent rank of tensor factors is the upper bound of the latent rank of the original tensor. More specifically, our tensor ring overlapped low-rank factor (TR-OLRF) model is formulated as follows:

\begin{equation}
\label{o_mod}
\begin{aligned}
\min \limits_{\{\tensor{G}_n\}_{n=1}^N,\tensor{X}} \ \sum_{n=1}^N\sum_{i=3}^N \Vert \Gamma_i(\mat{G}_n) \Vert_* + \frac{\lambda}{2}\Vert \tensor{X}- \tensor{Z}(\{\tensor{G}_n \}_{n=1}^N)\Vert_F^2 , \ s.t. \ P_\Omega(\tensor{X})=P_\Omega(\tensor{T}).
\end{aligned}
\end{equation}

The TR latent low-rank factor (TR-LLRF) model is outlined below:
\begin{equation}
\label{l_mod}
\begin{aligned}
\min \limits_{\{\tensor{G}_n\}_{n=1}^N, \tensor{X}}\ & \sum_{n=1}^N \sum_{i=1}^3 \Vert \Gamma_i(\mat{W}_{ni}) \Vert_* + \frac{\lambda}{2}\Vert \tensor{X}- \tensor{Z}(\{\tensor{G}_n \}_{n=1}^N)\Vert_F^2 , \\ &s.t. \ P_\Omega(\tensor{X})=P_\Omega(\tensor{T}),\ \tensor{G}_n=\sum\nolimits_{i=1}^3 \tensor{W}_{ni},\ n=1,\ldots,N.
\end{aligned}
\end{equation}
The two models have two distinctive advantages. Firstly, the low-rank assumption is placed on tensor factors instead of on the original tensor, this reduces the computational complexity of the SVD operation largely. Secondly, low-rankness on tensor factors can enhance the robustness to rank selection.

\subsection{Solving scheme}
\subsubsection{TR-OLRF}
	
To solve the equations (\ref{o_mod}) and (\ref{l_mod}), we apply the augmented Lagrangian multiplier method (ADMM) which is efficient and widely used. Because the variables of TR-OLRF are interdependent, we adopt alternative variables, and the augmented Lagrangian function of TR-OLRF model is:
\begin{equation}
\label{o_Lfunc}
\begin{aligned}
&\tensor{L}	\left( \{\tensor{G}_n\}_{n=1}^N,\tensor{X},\{\tensor{M}_{ni}\}_{n=1,i=1}^{N,3},\{\tensor{Y}_{ni}\}_{n=1,i=1}^{N,3} \right)=\sum_{n=1}^N\sum_{i=1}^3\Vert \Gamma_i(\mat{M}_{ni}) \Vert_*\\
&+\frac{\lambda}{2}\Vert  \tensor{X}- \tensor{Z}(\{\tensor{G}_n \}_{n=1}^N)\Vert_F^2 +\sum_{n=1}^N\sum_{i=1}^3 (<\tensor{Y}_{ni}, \tensor{M}_{ni}-\tensor{G}_n>+\frac{\mu}{2}\Vert  \tensor{M}_{ni}-\tensor{G}_n\Vert_F^2 ),
\end{aligned}
\end{equation}
where $\{\tensor{M}_{ni}\}_{n=1,i=1}^{N,3}$ are the alternative variables of $\{\tensor{G}_n\}_{n=1}^N$, $\{ \tensor{Y}_{ni}\}_{n=1,i=1}^{N,3}$ are Lagrangian multipliers, $<\cdot,\cdot>$ denotes the inner product, and $\mu>0$ is a penalty parameter.

To update ${\tensor{G}_n}$ and ${\tensor{M}_{ni}}$, $n=1,\ldots,N$, the augmented Lagrangian function are formulated by:
\begin{equation}
\label{tr_lG_o}
\begin{aligned}
\tensor{L}(\tensor{G}_n)=\sum_{i=1}^3\frac{\mu}{2}\Vert  \tensor{M}_{ni}-\tensor{G}_n+\frac{\tensor{Y}_{ni}}{\mu}  \Vert^2_F+\frac{\lambda}{2}\Vert \tensor{X}-\tensor{Z}(\{\tensor{G}_n\}_{n=1}^N) \Vert^2_F,
 \end{aligned}
\end{equation}

\begin{equation}
\label{tr_lM_o}
\begin{aligned}
\tensor{L}(\tensor{M}_{ni})=\sum_{i=1}^3\Vert \Gamma_i(\mat{M}_{ni}) \Vert_*+\sum_{i=1}^3\frac{\mu}{2}\Vert \tensor{M}_{ni}-(\tensor{G}_n-\frac{\tensor{Y}_{ni}}{\mu}) \Vert_F^2, i=1,2,3.
\end{aligned}
\end{equation}

For $n=1, \ldots, N$, the $k$th iteration update scheme of alternating direction method of multipliers (ADMM) of TR-OLRF model is listed below:
{\setlength\abovedisplayskip{5pt}
\setlength\belowdisplayskip{1pt}
\begin{equation}  
\label{o_all_u}
\left\{  
             \begin{array}{lr}  
             \tensor{G}_n^{k+1}=\bar{\Gamma}_2((\lambda \Delta_n(\mat{X}^k)(\mat{Q}_n^k)^T+\mu \sum_{i=1}^{3}\Gamma_2(\mat{M}_{ni}^k) +\sum_{i=1}^3\Gamma_2(\mat{Y}_{ni}^k))(\lambda\mat{Q}_n^k(\mat{Q}_n^k)^T+3\mu I)^{-1})  \\  
             \tensor{M}_{ni}^{k+1}=\bar{\Gamma}_i(D_{\frac{1}{\mu}}(\Gamma_i({\mat{G}_n^k}-\frac{1}{\mu}\mat{Y}_{ni}^k))),\ i=1,2,3\\  
            \tensor{X}^{k+1}=P_{\Omega}(\tensor{T})+P_{\bar{\Omega}}(\tensor{Z}(\{\tensor{G}_n^k\}_{n=1}^N))  \\
            \tensor{Y}_{ni}^{k+1}=\tensor{Y}_{ni}^k+\mu(\tensor{M}_{ni}^k-\tensor{G}_n^k), \ i=1,2,3\\
            \mu^{k+1}=max(\rho\mu^k,\mu_{max}).
             \end{array}  
\right.  
\end{equation}  }
where $\mat{Q}_n^k=(\Delta_2(\mat{G}_{\neq n}^{k}))^T$, $\bar{\Gamma}_n(\cdot)$ is the reverse operator of $\Gamma_n(\cdot)$ that transforms mode-$n$ matricization of a tensor to the original tensor, $D_{\beta}(\cdot)$ is the singular value thresholding (SVT) operator, i.e., if $\mat{U}\mat{S}\mat{V}^T$ is the singular value decomposition of matrix $\mat{A}$, then $D_\beta(\mat{A})=\mat{U}max\{\mat{S}-\beta \mat{I},0\}\mat{V}^T$, and $\bar{\Omega}$ is the set of indices of missing entries.

\subsubsection{TR-LLRF}
Similarly, the augmented Lagrangian function of TR-LLRF model can be written as:
\begin{equation}
\label{o_Lfunc}
\begin{aligned}
&\tensor{L}\left(\{\tensor{G}_n\}_{n=1}^N,\tensor{X},\{\tensor{W}_{ni}\}_{n=1,i=1}^{N,3},\{\tensor{Y}_n\}_{n=1}^N\right)= \sum_{n=1}^N\sum_{i=1}^3\Vert \Gamma_i(\mat{W}_{ni}) \Vert_*\\
&+\frac{\lambda}{2}\Vert  \tensor{X}- \tensor{Z}(\{\tensor{G}_n \}_{n=1}^N)\Vert_F^2+\sum_{n=1}^N (<\tensor{Y}_n, \sum_{i=1}^3\tensor{W}_{ni}-\tensor{G}_n)>+\frac{\mu}{2}\Vert \sum_{i=1}^3\tensor{W}_{ni}- \tensor{G}_n\Vert_F^2 ),
\end{aligned}
\end{equation}
To update ${\tensor{G}_n}$ and ${\tensor{W}_{ni}}$,  the augmented Lagrangian function is formulated by:
\begin{equation}
\label{tr_lG_l}
\begin{aligned}
\tensor{L}(\tensor{G}_n)=\frac{\mu}{2}\Vert  \sum_{i=1}^3\tensor{W}_{ni}-\tensor{G}_n+\frac{\tensor{Y}_{n}}{\mu}  \Vert^2_F+\frac{\lambda}{2}\Vert \tensor{X}-\tensor{Z}(\{\tensor{G}_n\}_{n=1}^N) \Vert^2_F,
 \end{aligned}
\end{equation}

\begin{equation}
\label{tr_lM_l}
\begin{aligned}
\tensor{L}(\tensor{W}_{ni})=\Vert \Gamma_i(\mat{W}_{ni}) \Vert_*+\frac{\mu}{2}\Vert \tensor{W}_{ni}-(\tensor{G}_n-\frac{\tensor{Y}_{n}}{\mu}-\sum_{j=1, j \neq i}^3\tensor{W}_{ni}) \Vert_F^2, i=1,2,3.
\end{aligned}
\end{equation}
The corresponding update scheme of TR-LLRF model is listed below:
{\setlength\abovedisplayskip{5pt}
\setlength\belowdisplayskip{1pt}
\begin{equation}  
\label{l_all_u}
\left\{
             \begin{array}{lr}  
             \tensor{G}_n^{k+1}=\bar{\Gamma}_2((\lambda \Delta_n(\mat{X}^k)(\mat{Q}_n^k)^T+\mu\sum_{i=1}^3 \Gamma_2(\mat{W}_{ni}^k)+\Gamma_2(\mat{Y}_n^k))(\lambda\mat{Q}_n^k(\mat{Q}_n^k)^T+\mu I)^{-1}) \\  
             \tensor{W}_{ni}^{k+1}=\bar{\Gamma}_i(D_{\frac{1}{\mu}}(\Gamma_i({\mat{G}_n^k}+\frac{1}{\mu}\mat{Y}_n^k-\sum_{j=1,j\neq i}^3\mat{W}_{ni})), i=1,2,3\\  
            \tensor{X}^{k+1}=P_{\Omega}(\tensor{T})+P_{\bar{\Omega}}(\tensor{Z}(\{\tensor{G}_n^k\}_{n=1}^N))  \\
            \tensor{Y}_n^{k+1}=\tensor{Y}_n^k+\mu(\sum_{i=1}^3 \tensor{W}_{ni}^k-\mat{G}_n^k)\\
            \mu^{k+1}=max(\rho\mu^k,\mu_{max}),
             \end{array}  
\right.  
\end{equation}  }

The ADMM solving model is updated iteratively based on the above model and updating scheme. The implementation process and hyper-parameter selection of the two algorithms are summarized in Alg. 1 and Alg. 2.
\begin{table}[H]
\centering
\footnotesize
\resizebox{\textwidth}{25mm}{
\begin{tabular}{p{7cm}|p{7cm}}
\hline
\textbf{Alg. 1}  TR overlapped low-rank factors (TR-OLRF) & \textbf{Alg. 2} TR latent low-rank factors (TR-LLRF)\\
\hline
1:\;\; \textbf{Input}: $P_{\Omega}(\tensor{T})$, initial TR-rank $\vect{r}$ .  & 1:\;\;\textbf{Input}: $P_{\Omega}(\tensor{T})$, initial TR-rank $\vect{r}$,.\\
2:\;\; \textbf{Initialization}: $k=0$, $\mu^0=1$, $\mu_{max}=10^2$, $\rho=1.01$, $tol=10^{-6}$, $\lambda=10$, element of $\tensor{G}_n^0$ s.t. $\tensor{N}\sim(0,1)$, $\tensor{Y}_{ni}^0=0$, $\tensor{M}_{ni}^0=0$, $n=1,\ldots,N$, $i=1,2,3$.  & 2:\;\; \textbf{Initialization}: $k=0$, $\mu^0=1$, $\mu_{max}=10^2$, $\rho=1.01$, $tol=10^{-6}$, $\lambda=10$, element of $\tensor{G}_n^0$ s.t. $\tensor{N}\sim(0,1)$, $\tensor{Y}_{n}^0=0$, $\tensor{W}_{ni}^0=0$, $n=1,\ldots,N$, $i=1,2,3$. \\
3:\;\; \textbf{While} the stopping condition is not satisfied \textbf{do}&3:\;\; \textbf{While} the stopping condition is not satisfied \textbf{do}\\
4:\;\;\; k=k+1;&4:\;\;\; k=k+1;\\
5:\;\;\; Update variables by equation (\ref{o_all_u}).&5:\;\;\; Update variables by equation (\ref{l_all_u}).\\
6:\;\;\; \textbf{If} $\Vert \tensor{X}^{k+1}-\tensor{X}^k \Vert_F / \Vert  P_{\Omega}(\tensor{T})\Vert_F<tol$, \textbf{break}&6:\;\;\; \textbf{If} $\Vert \tensor{X}^{k+1}-\tensor{X}^k \Vert_F / \Vert  P_{\Omega}(\tensor{T})\Vert_F<tol$, \textbf{break}\\
7: \;\;\;\textbf{End while} &7: \;\;\;\textbf{End while} \\
8: \;\;\;\textbf{Output}: $\tensor{X}$ and $\tensor{G}_n$, $n=1,\ldots,N$.&8: \;\;\;\textbf{Output}: $\tensor{X}$ and $\tensor{G}_n$, $n=1,\ldots,N$.\\
\hline
\end{tabular}
}
\end{table}

\subsection{Computational complexity}

\begin{table}[H]
\caption{Computational complexity}
\label{acc}
\footnotesize
\begin{center}
\begin{tabular}{c|c}
\hline
Algorithm &Computational Complexity \\
\hline
TR-OLRF&$\tensor{O}(NI^3R^3+NIR^4+NR^6)$ \\
TR-LLRF&$\tensor{O}(NI^3R^3+NIR^4+NR^6)$ \\
 TR-ALS&$\tensor{O}(NMR^4+NR^6)$\\
 TT-SiLRTC&$\tensor{O}(NI^{3/2N}+NI^{3/2N-1})$ \\
 SiLRTC &$\tensor{O}(NI^{N+1})$\\
 BCPF&$\tensor{O}(NMR^2+R^3)$\\
\hline
\end{tabular}
\end{center}
\end{table}

We next compared the computational complexity of our TR-OLRF and TR-LLRF to the state-of-the-art algorithms TR-ALS \cite{wang2017efficient}, SiLRTC-TT \cite{bengua2017efficient}, SiLRTC \cite{liu2013tensor} and FBCP \cite{zhao2015bayesian}. The comparative algorithms are state-of-the-art algorithms and are similar to our algorithms. The complexities are summarized in Tab. \ref{acc}, where we denote the dimension of tensor by $N$, $I_1=\cdots=I_N=I$, and all the TT-ranks, TR-ranks and CP ranks are set to $R$. From Tab. \ref{acc} we can see that compared to Schatten norm based algorithms, the computational complexity of our algorithms are linear in tensor dimension. Compared to TR-ALS and BCPF, the complexity of our algorithms is independent from the number of observed entries. The computational complexity of our algorithms increase fast when $R$ increases, however, due to the linear scalability of TR decomposition, $R$ is often small in model selection of proposed algorithms. In addition, most of the stated algorithms are rank adaptive, i.e., robust to rank selection.

\section{Experiment results}

\subsection{Synthetic data}
To verifying the performance of our two proposed algorithms, we test two tensors of size $10\times10\times10\times10$ and $4 \times4\times4\times6\times6\times6$. The tensors were generated by TR factors of TR-ranks $\{4,5,4,5\}$ and $\{4,4,4,4,4,4\}$ respectively. The values of the TR factors were drawn from an $i.i.d.$ normal distribution $\tensor{N} \sim (0,0.5)$. We define $SSR$ as the sum of square root of TR-rank (i.e. $SSR=\sum_{n=1}^N\sqrt{R_n}$) to be the index of model complexity. The observed entries of the tensors were randomly removed. We verified the performance of the proposed two algorithms in several scenarios, with the mean RSE values of 10 times of dependent experiments as the final results. All the hyper-parameters of the two algorithms were set according to Alg. 1 and Alg. 2. 

\begin{figure}[htb]
\begin{center}
\includegraphics[width=1\linewidth]{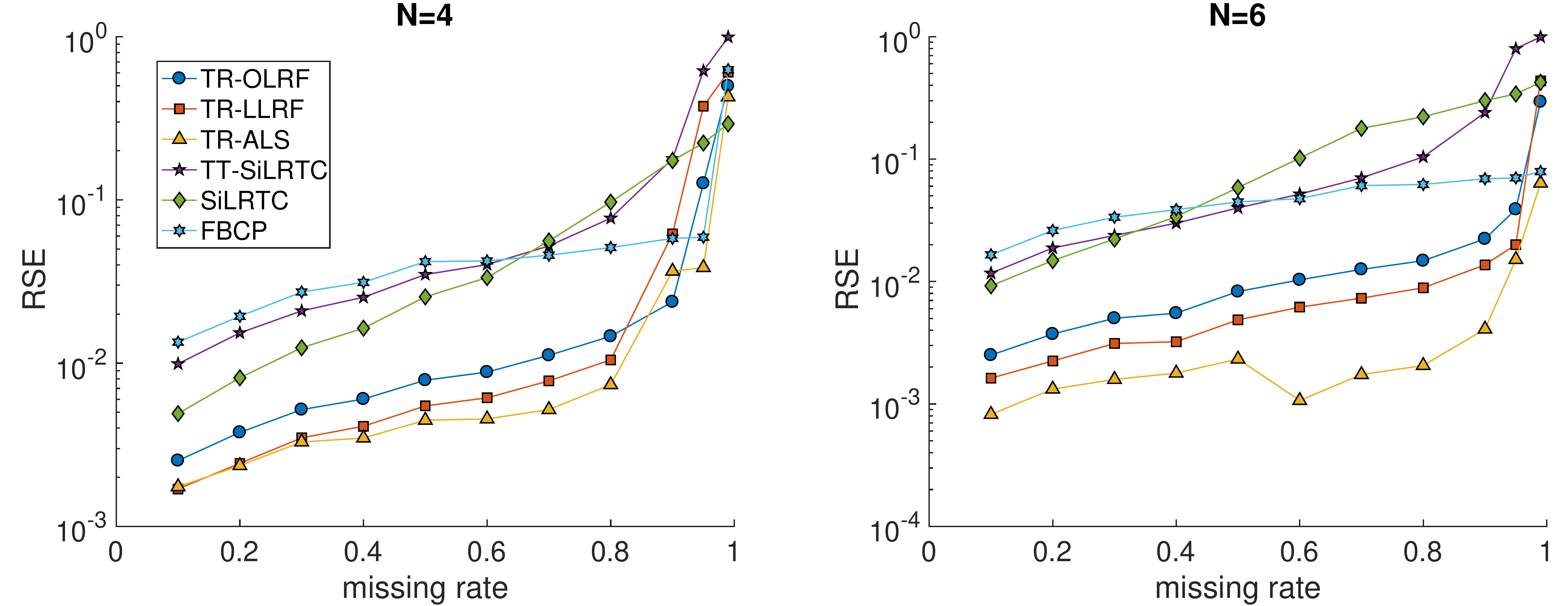}
\label{sim}
\end{center}
\caption{Completion performance of six algorithms under different missing rates}\label{simu_mr}
\end{figure}

For the first experiment, we test the completion performance of our two algorithms and four other state-of-the-art algorithms under different missing rates, from 0.1 to 0.99. For our algorithms, we set the TR-rank to be the same as the real rank of the synthetic data and other hyper-parameters were set as default. For other compared algorithms, we tuned the hyper-parameters respectively to obtain the best results of each algorithm. Fig. \ref{simu_mr} shows the experiment results for the order-four tensor and order-six tensor respectively. 

\begin{figure}[htb]
\label{simu_3}
\begin{center}
\includegraphics[width=1\linewidth]{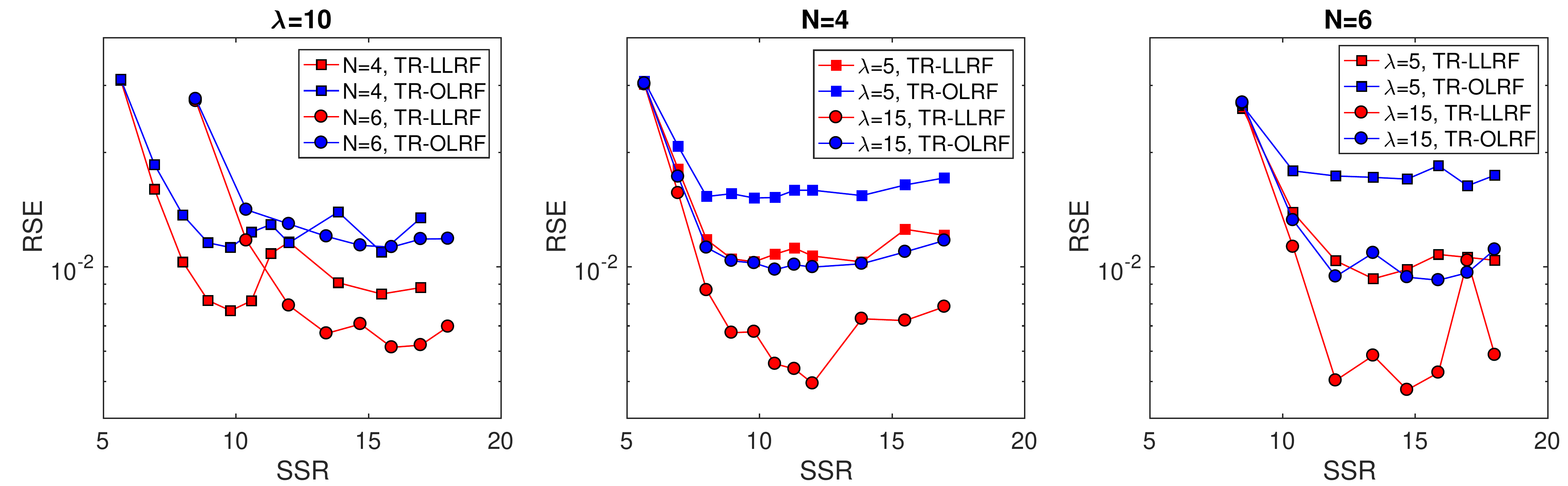}
\end{center}
\caption{Algorithm robustness to rank and $\lambda$}\label{sim}
\end{figure}

For the second experiment, we tested the completion performance of our two algorithms under various SSR values, the missing rate was set to $0.7$, and used again two different tensors. The results in the first picture of Fig. \ref{simu_3} show that our two algorithms obtained the lowest RSE values when the SSR was near the real SSR, and when the SSR value increased, the RSE value remained stable. This indicates that our algorithms are robust to rank selection.

For the third and forth experiments, we tested the performance of our algorithms over different values of $\lambda$, missing rate was set to $0.7$, and TR-rank is chosen as the real rank of the two tensors. Fig. \ref{simu_3} shows the robustness for the three different values of $\lambda$ and verifies that our two algorithms are robust to the selection of $\lambda$.

\subsection{Hyperspectral image}
A hyperspectral image of size $200\times\ 200\times80$ was next considered. This was an image of urban landscape collected by a satellite. We compare our TR-OLRF and TR-LLRF to TR-ALS , TT-SiLRTC , SiLRTC and BCPF. We examined order-three, order-five, order-seven and order-eight tensors respectively. The missing rate is set as 0.9 and the hyper-parameters are set as defaults. For each tensor, we choose all the TR-ranks as a same value, i.e., $R_1=\ldots=R_N$. The tensor size and TR-ranks are recorded in the first column of Tab. \ref{hsi}, and the RSE values of each tensor against each algorithm are listed in Tab. \ref{hsi}.

From the results we can see, our algorithms significantly outperform TT-SiLRTC, SiLRTC, BCPF. Though the results of TR-ALS are comparable to our algorithms, it should be noted that the 
computational time of TR-ALS is more than double of the time TR-OLRF and TR-LLRF spent (1891 seconds vs 756 seconds and 988 seconds) in order to get the similar results.
\begin{table}[H]
\caption{Completion results under four different tensor dimensions}
\label{hsi}
\footnotesize
\resizebox{\textwidth}{10mm}{
\begin{tabular}{c|c|c|c|c|c|c}
\hline
 & TR-OLRF & TR-LLRF & TR-ALS & TT-SiLRTC & SiLRTC &BCPF\\
\hline
$\mathbb{R}^{200\times\ 200\times80},R_n=12$&0.0710 &\textbf{0.0677} &0.0681  &0.4572 &0.3835&0.3750\\
$\mathbb{R}^{10\times 20\times 10 \times 20 \times 80}, R_n=18$&\textbf{0.1062} &0.1072 &0.1122 &0.4895 &0.4307 &0.3742\\
$\mathbb{R}^{4\times5\times 10 \times 4 \times5 \times 10 \times 80}, R_n=20$&\textbf{0.1436} &0.1483 &0.1497 &0.5051 &0.4408 &0.3680\\
$\mathbb{R}^{8\times 5\times 5 \times 8 \times 5\times 5\times 8\times 10}, R_n=22$& 0.1520&0.1524 &\textbf{0.1430} &0.4957&0.4526 &0.3981\\
\hline
\end{tabular}
}
\end{table}

\section{Conclusion}
In order to solve the large-scale SVD calculation and rank selection problem that most tensor completion methods have. We proposed two algorithms which impose low-rank assumption on tensor factors. Based on tensor ring decomposition, we proposed two optimization models named as TR-OLRF and TR-LLRF. The two models can be solved efficiently by ADMM algorithm. We test the algorithms on synthetic data in various situations by synthetic data and real world data. The high performance and high efficiency of ur algorithms are obtained from the experiment results. In addition, the results also show that the proposed algorithms are robust to tensor rank and other model parameters. The proposed method is heuristic to all the model-based low-rank tensor completion and decomposition, and it can be applied to various tensor decompositions to create more efficient and robust algorithms. 

\section*{Acknowledgement}
\bibliographystyle{plain}
\bibliography{paper_arxiv}

\end{document}